\documentclass[12pt,a4paper]{article}

\usepackage{graphicx}
\usepackage{amsmath,amsfonts}
\usepackage{amsthm} % theorem without numbering
\usepackage{mathrsfs}
\usepackage[margin=0.75in]{geometry}

\usepackage{color} 
\usepackage{hyperref}
\hypersetup{
    colorlinks,
    citecolor=green,
    filecolor=blue,
    linkcolor=red,
    urlcolor=blue
}

\theoremstyle{plain} % theorem without numbering

\newtheorem*{theorem*}{Theorem} % theorem without numbering
\newtheorem*{proposition*}{Proposition} % theorem without numbering
\newtheorem*{lemma*}{Lemma}

\newtheorem*{definition*}{Definition}

\def\Pb{\mathbf{P}}
\def\Ex{\mathbf{E}}

\def\d{{\rm d}} %S
\def\F{{\mathcal F}} %S
\def\R{ {\rm {\bf R}}} %S
\def\X{ {\rm {\bf X}}} %S
\def\N{{\mathcal N}} %S
\def\C{{\mathcal C}} %S
\def\Z{{\mathcal Z}} %S
\def\L{{\rm L}} %S
\begin{document}

\title{{\bf Second Order Asymptotic Efficiency for a Poisson Process\thanks{Completely revised version of the originally published paper}}}

\author{Samvel B. GASPARYAN\thanks{\href{mailto:GasparyanSB@gmail.com}{GasparyanSB@gmail.com}}\\
\small{Universit\'e du Maine, Le Mans, France}\\
\small{Yerevan State University, Yerevan, Armenia}\\
}
\date{}

\maketitle
\begin{abstract}
\begin{center} 
%{\bf \formatdate{04}{02}{2015}} 
\end{center}
We consider the problem of the estimation of the mean function of an inhomogeneous Poisson process when its intensity function is periodic. For the mean integrated squared error (MISE) there is a classical lower bound for all estimators and the empirical mean function attains that lower bound, thus it is asymptotically efficient. Following the ideas of the work by Golubev and Levit, we compare asymptotically efficient estimators and propose an estimator which is second order asymptotically efficient. Second order efficiency is done over Sobolev ellipsoids, following the idea of Pinsker.

\end{abstract}

\textbf{MSC2010 numbers:} 62G05, 62M05 

\textbf{Keywords:} Poisson process, second order estimation, asymptotic efficiency.

%\tableofcontents

\section{Introduction}
We consider the problem of non-parametric estimation of the mean function of an inhomogeneous Poisson process. We suppose that the unknown intensity function is periodic. It is known that the empirical mean function is an asymptotically efficient (in several senses, see e.g. Kutoyants \cite{K1},\cite{K2}) estimator. Particularly, we are interested in asymptotic efficiency with respect to the mean integrated squared error (MISE). Note that there are many estimators that are asymptotically efficient in this sense. The goal of the present work is to choose in this class of asymptotically efficient estimators an estimator which is asymptotically efficient of the second order. Such a statement of the problem was considered by Golubev and Levit \cite{GL} in the problem of distribution function estimation for the model of independent, identically distributed random variables. Then the ideas of this work were applied to the second order asymptotically efficient estimation for different models - Dalalyan, Kutoyants \cite{DK} proved second order asymptotic efficiency in the estimation problem of the invariant density of an ergodic diffusion process, in partial linear models the second order asymptotic efficiency was proved by Golubev, H\"ardle \cite{GH}. In this paper (combined with the paper \cite{G}) we prove the second order asymptotic efficiency result for the mean function of a Poisson process. The main idea that led to development of these type of problems was proposed by Pinsker in \cite{P} (more details on the Pinsker bound can be found in \cite{N}, \cite{T}).

\section{Auxiliary Results}

We are given a probability space $(\Omega,\F,\Pb)$ and a stochastic process $\X^T=\{X_t,\,t\in[0,T]\}.$ Recall that $\X^T$ is an inhomogeneous Poisson process if $1.\,X_0=0$ a.s. $2.$ The increments of the process $\X^T$ on the disjoints intervals are independent random variables. 3. We have
$$ \Pb(X_t-X_s=k)=\frac{\left[\Lambda(t)-\Lambda(s)\right]^k}{k!}e^{-\left[\Lambda(t)-\Lambda(s)\right]},\, 0\leq s<t\leq T,k\in {\mathcal Z}_+.$$
Here $\{\Lambda(t),\,t\in[0,T]\}$ is a non-decreasing function, and is called the mean function of the Poisson process, because $\Ex X(t)=\Lambda(t)$. If the mean function is absolutely continuous 
$$\Lambda(t)=\int_0^t \lambda(s)\d s,$$
then $\{\lambda(t),\, 0\leq t\leq T\}$ is called the intensity function.

Let us consider the problem of estimation $\Lambda(\cdot),$ when its intensity function is a $\tau$-periodic function. For simplicity we suppose that $T=T_n=\tau n.$ Then the observations $\X^T=\{X_t\,,t\in[0,\tau n]\},$ can be written in the form
\begin{equation}
\label{eq:1}
\X^n=(\X_1,\X_2,\cdots,\X_n),
\end{equation}
where
$$\X_j=\{X_j(t),\, 0\leq t\leq\tau\},\ \ X_j(t)=X_{(j-1)\tau+t}-X_{(j-1)\tau},\,j=1,\cdots,n.$$
It is well known that \textit{the empirical estimator}
$$\hat\Lambda_n(t)=\frac{1}{n}\sum_{j=1}^n X_j(t),\, t\in[0,\tau]$$
is consistent and asymptotically normal: for all $t\in[0,\tau]$
$$\sqrt{n}(\hat\Lambda_n(t)-\Lambda(t))\Longrightarrow\N(0,\Lambda(t)).$$
Moreover, this estimator is asymptotically efficient in the sense of the following lower bound: for all estimators $\{\bar\Lambda_n(t),\, t\in[0,\tau]\}$ and all $t^\ast\in[0,\tau]$ we have

$$\lim_{\delta\rightarrow 0}\varliminf_{n\rightarrow +\infty}\sup_{\Lambda\in V_\delta}n \Ex_\Lambda(\bar\Lambda_n(t^\ast)-\Lambda(t^\ast))^2\geq \Lambda^\ast(t^\ast),$$
where $V_\delta=\{\Lambda(\cdot):\ \ \sup_{t\in[0,\tau]}|\Lambda(t)-\Lambda^\ast(t)|\leq \delta\}$ and for the empirical mean function
one has an equality. This is a particular case of a general lower bound given in Kutoyants \cite{K1}. Similar inequality holds for MISE (\cite{K2})

\begin{equation}
\label{eq:2}
\lim_{\delta\rightarrow 0}\varliminf_{n\rightarrow +\infty}\sup_{\Lambda\in V_\delta}n \int_0^\tau \Ex_\Lambda(\bar\Lambda_n(s)-\Lambda(s))^2\d s\geq\int_0^\tau\Lambda^\ast(s)\d s.
\end{equation}

\begin{definition*}
The estimators $\Lambda_n^\ast(\cdot)$ for which we have an equality in \eqref{eq:2}, i.e., 
$$\lim_{\delta\rightarrow 0}\lim_{n\rightarrow +\infty}\sup_{\Lambda\in V_\delta}n \int_0^\tau \Ex_\Lambda(\Lambda_n^\ast(s)-\Lambda(s))^2\d s=\int_0^\tau\Lambda^\ast(s)\d s,
$$
are called (first order) asymptotically efficient.
\end{definition*}
The empirical mean function is an asymptotically efficient estimator also in this sense. (\cite{K2})

The goal of the present work is to find in the class of first order asymptotically efficient estimators an estimator which is second order asymptotically efficient. We follow the mains steps of the proof of Golubev, Levit \cite{GL}.

\section{Main Result}

For a given integer $m>1$ consider the following set of non-decreasing, positive functions on $[0,\tau]$ such that their $(m-1)$th derivative is absolutely continuous and
\begin{align}
\label{15}
\F_m^{per}(R,S)=\left\{ \Lambda(\cdot):\, \int_0^\tau[\Lambda^{(m)}(t)]^2 \d t\leq R,\,\Lambda(0)=0,\,\Lambda(\tau)=S\right\},\,m>1,
\end{align}
where $R>0,\,S>0$ are given constants. Periodicity of the Poisson process means that the intensity function $\lambda(\cdot)$ is periodic, hence the equality of its values and the values of its derivatives on the endpoints of the interval $[0,\tau]$ (for estimating a non-periodic function, see, for example, \cite{DH}). Introduce as well 
\begin{align}
\label{19}
\Pi=\Pi_m(R,S)=(2m-1)R\left(\frac{S}{\pi R}\frac{m}{(2m-1)(m-1)}\right)^\frac{2m}{2m-1}.
\end{align}
\begin{proposition*}
\label{P3} Consider Poisson observations $\X=(\X_1,\X_2,\cdots,\X_n)$ defined in \eqref{eq:1}. Then, for all estimators $\bar\Lambda_n(t)$ of the mean function $\Lambda(t),$ following lower bound holds
$$\varliminf_{n\rightarrow+\infty}\sup_{\Lambda\in {\mathcal F}_m(R,S)}n^{\frac{2m}{2m-1}}\left(\int_0^\tau \Ex_\Lambda  (\bar\Lambda_n(t)-\Lambda(t))^2\d t-\frac{1}{n}\int_0^\tau \Lambda(t)\d t\right)\geq-\Pi.$$
\end{proposition*}
This proposition is going to be presented in the forthcoming work \cite{G} (proof relies on a method developed in \cite{GiL}). In the next theorem we propose an estimator which attains this lower bound, thus we prove that this lower bound is sharp. Introduce
\begin{align*}
\Lambda_n^\ast(t)=\hat\Lambda_{0,n}\phi_0(t)+\sum_{l=1}^{N_n}\tilde K_{l,n}\hat\Lambda_{l,n}\phi_{l}(t),
\end{align*}
where $\{\phi_l\}_{l=0}^{+\infty}$ is the trigonometric cosine basis in $\L_2[0,\tau]$ (see \eqref{30} below), $\hat\Lambda_{l,n}$ are the Fourier coefficients of the empirical mean function with respect to this basis and 
\begin{align*}
&\tilde K_{l,n}=\left(1-\left|\frac{\pi l}{\tau}\right|^m \alpha_n^\ast\right)_+,\ \ \alpha_n^\ast=\left[\frac{S}{n R}\frac{\tau}{\pi}\frac{m}{(2m-1)(m-1)}\right]^\frac{m}{2m-1},\\
&N_n=\frac{\tau}{\pi}(\alpha_n^\ast)^{-\frac{1}{m}}\approx {\rm C} n^\frac{1}{2m-1}, \ \ x_+=\max(x,0),\, x\in\R.
\end{align*}
The main result of this paper states
\begin{theorem*} 
\label{P4}
The estimator $\Lambda_n^\ast(t)$ attains the lower bound described above, that is, 
$$\lim_{n\rightarrow+\infty}\sup_{\Lambda\in {\mathcal F}_m(R,S)}n^{\frac{2m}{2m-1}}\left(\int_0^\tau \Ex_\Lambda  (\Lambda_n^\ast(t)-\Lambda(t))^2\d t-\frac{1}{n}\int_0^\tau \Lambda(t)\d t\right)=-\Pi.$$
\end{theorem*}

\section{The Proof}
Consider the $\L_2[0,\tau]$ Hilbert space. Evidently, $\F_m^{per}(R,S)\subset\L_2[0,\tau].$ The main idea of the proof is to replace the estimation problem of the infinite-dimensional (continuum) mean function by the estimation problem of infinite-dimensional but countable vector of its Fourier coefficients. Recall that the space $\L_2[0,\tau]$ is isomorphic to the space
$$\ell_2=\left\{\theta=(\theta_l)_{l=0}^{+\infty}:\quad \sum_{l=0}^{+\infty}\theta_l^2<+\infty\right\},\ \ ||\theta||=\left(\sum_{l=0}^{+\infty}\theta_l^2\right)^\frac{1}{2}.$$
Our first goal is to describe the set $\Theta\subset\ell_2$ of Fourier coefficients of the functions from the set $\F_m^{per}(R,S).$ 

Consider a complete, orthonormal system in the space $\L_2[0,\tau],$
\begin{align}
\label{30}
\phi_0(t)=\sqrt{\frac{1}{\tau}},\, \phi_{l}(t)=\sqrt{\frac{2}{\tau}}\cos\left(\frac{\pi l}{\tau}t\right),\,l\in\N.
\end{align}
Each function $f\in\L_2[0,\tau]$ is a $\L_2-$limit of its Fourier series
$$f(t)=\sum_{l=0}^{+\infty}\theta_l\phi_l(t),\quad \theta_l=\int_0^\tau f(t)\phi_l(t)\d t.$$
Suppose that 
$$\Lambda_l=\int_0^\tau \Lambda(t)\phi_l(t)\d t,\quad \lambda_l=\int_0^\tau \lambda(t)\phi_l(t)\d t.$$
Then
\begin{lemma*}
\label{L1} 
The mean function $\Lambda$ belongs to the set $\F_m^{per}(R,S)$ (see \eqref{15}) if and only if its Fourier coefficients w.r.t. the cosine trigonometric basis satisfy
\begin{align}
\label{35}
\sum_{l=1}^{+\infty}\left(\frac{\pi l}{\tau}\right)^{2m}\Lambda_l^2\leq R,\ \ \Lambda(\tau)=S,
\end{align}
or, the Fourier coefficients of its intensity function satisfy
\begin{align}
\label{34}
\sum_{l=1}^{+\infty}\left(\frac{\pi l}{\tau}\right)^{2(m-1)}\lambda_l^2\leq R,\ \ \Lambda(\tau)=S.
\end{align}

\end{lemma*}
For the proof see, for example, \cite{T} (Lemma A.3). 
To introduce the estimator denote the Fourier coefficients of the empirical mean function by
$$\hat\Lambda_{l,n}=\int_0^\tau\hat\Lambda_n(t)\phi_l(t)\d t,\,l\in\Z_+,\ \ \hat\Lambda_n(t)=\frac{1}{n}\sum_{j=1}^n X_j(t).$$
Consider the estimator 
$$\tilde\Lambda_n(t)=\sum_{l=0}^{+\infty}\tilde\Lambda_{l,n}\phi_l(t),\,\tilde\Lambda_{l,n}=K_{l,n}\hat\Lambda_{l,n}.$$
Here $K_{l,n}$ are some numbers. Without loss of generality we can take $K_{0,n}=1,$ that is
$\tilde\Lambda_{0,n}=\hat\Lambda_{0,n}.$ In this case, using the Parseval's equality, we get

\begin{align}
\label{22}
&\Ex_\Lambda\|\tilde\Lambda_n-\Lambda\|^2-\Ex_\Lambda\|\hat\Lambda_n-\Lambda\|^2=\sum_{l=1}^{+\infty}(K_{l,n}^2-1)\sigma_{l,n}^2+\sum_{l=1}^{+\infty}\left|K_{l,n}-1\right|^2\Lambda_{l}^2.
\end{align}
Here $\sigma_{l,n}^2=\Ex_\Lambda(\hat\Lambda_{l,n}-\Lambda_{l})^2.$ To compute this quantity, introduce the notation
$$
\pi_j(t)=X_j(t)-\Lambda(t).
$$ 
In the sequel, we are going to use the following property of stochastic integrals (see, for example, \cite{K1},\cite{K2})
\begin{align*}
\Ex_\Lambda\bigg[\int_0^\tau f(t)\d\pi_j(t)\int_0^\tau g(t)\d\pi_j(t)\bigg]=\int_0^\tau f(t)g(t)\d\Lambda(t),\ \ f,g\in \L_2[0,\tau].
\end{align*}
Further, in view of the integration by parts, we have
\begin{align*}
\hat\Lambda_{l,n}-\Lambda_l = \frac{1}{n}\sum_{j=1}^n \int_0^\tau \pi_j(t)\phi_l(t)\d t= \frac{1}{n}\sum_{j=1}^n \int_0^\tau \left(\int_t^\tau\phi_l(s)\d s\right)\d\pi_j(t),
\end{align*}
which entails that 
\begin{align*}
\sigma_{l,n}^2=\Ex_\Lambda|\hat\Lambda_{l,n}-\Lambda_l|^2
    & = \frac{1}{n}\int_0^\tau \left(\int_t^\tau\phi_l(s)\d s\right)^2\d\Lambda(t).
\end{align*}
Simple algebra yields
\begin{align*}
\sigma_{l,n}^2=\frac{1}{n}\left(\frac{\tau}{\pi l}\right)^2\left[\Lambda(\tau)-\frac{2}{\tau}\int_0^\tau\cos\left(\frac{2\pi l}{\tau}t\right)\lambda(t)\d t\right].
\end{align*}
Combining with (\ref{22}), this leads to 
\begin{align}
\label{23}
&\Ex_\Lambda\|\tilde\Lambda_n-\Lambda\|^2-\Ex_\Lambda\|\hat\Lambda_n-\Lambda\|^2=\frac{S}{n}\sum_{l=1}^{+\infty}\left(\frac{\tau}{\pi l}\right)^2(K_{l,n}^2-1)\nonumber\\
&+\sum_{l=1}^{+\infty}\left(K_{l,n}-1\right)^2\Lambda_{l}^2+\frac{1}{n}\sqrt{\frac{2}{\tau}}\sum_{l=1}^{+\infty}\left(\frac{\tau}{\pi l}\right)^2(1-K_{l,n}^2)\lambda_{2l}.
\end{align}
For the third term in the right-hand side we have
\begin{align*}
&\left|\frac{1}{n}\sqrt{\frac{2}{\tau}}\sum_{l=1}^{+\infty}\left(\frac{\tau}{\pi l}\right)^2(1-K_{l,n}^2)\lambda_{2l}\right|\leq\\
&\leq\frac{1}{n}\sqrt{\frac{2}{\tau}}\max_l\frac{|1-K_{l,n}^2|}{\left(\frac{\pi l}{\tau}\right)^m}\sum_{l=1}^{+\infty}\left(\frac{\pi l}{\tau}\right)^{m-1}\lambda_{2l}\left(\frac{\pi l}{\tau}\right)^{-1}\\
&\leq\frac{1}{n}\sqrt{\frac{2}{\tau}}\max_l\frac{|1-K_{l,n}^2|}{\left(\frac{\pi l}{\tau}\right)^m}\left(\sum_{l=1}^{+\infty}\left(\frac{\pi l}{\tau}\right)^{2(m-1)}\lambda_{2l}^2\right)^\frac{1}{2}\left(\sum_{l=1}^{+\infty}\left(\frac{\pi l}{\tau}\right)^{-2}\right)^\frac{1}{2}.
\end{align*}
Using \eqref{34} from the Lemma  we obtain
$$\left(\sum_{l=1}^{+\infty}\left(\frac{\pi l}{\tau}\right)^{2(m-1)}\lambda_{2l}^2\right)^\frac{1}{2}\leq \sqrt{R}.$$
Hence
$$\left|\frac{1}{n}\sqrt{\frac{2}{\tau}}\sum_{l=1}^{+\infty}\left(\frac{\tau}{\pi l}\right)^2(1-K_{l,n}^2)\lambda_{2l}\right|
\leq\frac{{\rm C}}{n}\max_l\frac{|1-|K_{l,n}|^2|}{\left(\frac{\pi l}{\tau}\right)^m}$$
Now, consider the first two terms of the right-hand side of the equation \eqref{23}. Introduce a set of possible kernels (for all $c_n>0$)
$${\mathcal C_n}=\left\{K_{l,n}:|K_{l,n}-1|\leq \left|\frac{\pi l}{\tau}\right|^m c_n\right\}.$$
From \eqref{35} follows
\begin{align*}
\frac{S}{n}\sum_{l=1}^{+\infty}\left(\frac{\tau}{\pi l}\right)^2(K_{l,n}^2-1)+\sum_{l=1}^{+\infty}\left|K_{l,n}-1\right|^2\Lambda_{l}^2=&\frac{S}{n}\sum_{l=1}^{+\infty}\left(\frac{\tau}{\pi l}\right)^2(K_{l,n}^2-1)+\\
\sum_{l=1}^{+\infty}\frac{\left |K_{l,n}-1\right|^2}{\left(\frac{\pi l}{\tau}\right)^{2m}}\left(\frac{\pi l}{\tau}\right)^{2m}\Lambda_{l}^2\leq&\frac{S}{n}\sum_{l=1}^{+\infty}\left(\frac{\tau}{\pi l}\right)^2(K_{l,n}^2-1)+c_n^2R.
\end{align*}
Hence, minimizing the later over the set $\C_n$ 
\begin{equation}
\label{29}
\tilde K_{l,n}=\arg\min_{{\mathcal C}_n}|K_{l,n}|=\left(1-\left|\frac{\pi l}{\tau}\right|^m c_n\right)_+,
\end{equation}
we obtain
\begin{align}
\label{28}
&\sup_{\Lambda\in\F_m^{per}(R,S)}\left(\Ex_\Lambda\|\tilde\Lambda_n-\Lambda\|^2-\Ex_\Lambda\|\hat\Lambda_n-\Lambda\|^2\right)\leq\nonumber\\
&\frac{S}{n}\sum_{l=1}^{+\infty}\left(\frac{\tau}{\pi l}\right)^2(\tilde K_{l,n}^2-1)+c_n^2R+\frac{{\rm C}}{n}\max_l\frac{|1-\tilde K_{l,n}^2|}{\left(\frac{\pi l}{\tau}\right)^m}.
\end{align}
Here $\tilde\Lambda_n(t)$ is the estimator corresponding to the kernel $\tilde K(u).$ In fact, we have not yet constructed the estimator. We have to specify the sequence of positive numbers $c_n$ in the definition \eqref{29}. Consider the function
$$H(c_n)=\frac{S}{n}\sum_{l=1}^{+\infty}\left(\frac{\tau}{\pi l}\right)^2(\tilde K_{l,n}^2-1)+c_n^2R$$
and minimize it with respect to the positive sequence $c_n$. Introduce as well $N_n=\frac{\tau}{\pi}c_n^{-\frac{1}{m}}.$ Then
\begin{align*}
H(c_n)=\frac{S}{n}\left[\sum_{l\leq N_n}\left(\frac{\tau}{\pi l}\right)^2\left(c_n^2\left(\frac{\pi l}{\tau}\right)^{2m}-2c_n\left(\frac{\pi l}{\tau}\right)^{m}\right)-\sum_{l>N_n}\left(\frac{\tau}{\pi l}\right)^2\right]+c_n^2R.
\end{align*}
To minimize this function consider its derivative
\begin{align}
\label{24}
H'(c_n)=\frac{S}{n}\sum_{l\leq N_n}\left(\frac{\tau}{\pi l}\right)^2\left[2c_n\left(\frac{\pi l}{\tau}\right)^{2m}-2\left(\frac{\pi l}{\tau}\right)^{m}\right]+2c_n R=0.
\end{align}
Consider such sums $(\beta\in\N)$
\begin{align*}
\sum_{l\leq N_n}l^\beta=\sum_{l=1}^{[N_n]}\left(\frac{l}{[N_n]}\right)^\beta[N_n]^\beta=[N_n]^{\beta+1}\sum_{l=1}^{[N_n]}\left(\frac{l}{[N_n]}\right)^\beta\frac{1}{[N_n]},
\end{align*}
hence, if $c_n\longrightarrow0,$ as $n\longrightarrow+\infty,$
$$\frac{1}{[N_n]^{\beta+1}}\sum_{l\leq N_n}l^\beta\longrightarrow\int_0^1 x^\beta\d x,$$
that is,
$$\sum_{l\leq N_n}l^\beta=\frac{[N_n]^{\beta+1}}{\beta+1}(1+o(1)),\ \ n\longrightarrow+\infty.$$
Using this identity we can transform \eqref{28} (remembering that $N_n=\frac{\tau}{\pi}c_n^{-\frac{1}{m}}$)
\begin{align*}
&\frac{S}{n}\left(c_n\left(\frac{\pi}{\tau}\right)^{2(m-1)}\sum_{l\leq N_n}l^{2(m-1)}-\left(\frac{\pi}{\tau}\right)^{m-2}\sum_{l\leq N_n}l^{m-2}\right)=-c_n R,\\
&\frac{S}{n}\left(c_n\left(\frac{\pi}{\tau}\right)^{2(m-1)}\frac{N_n^{2m-1}}{2m-1}-\left(\frac{\pi}{\tau}\right)^{m-2}\frac{N_n^{m-1}}{m-1}\right)=-c_n R(1+o(1)),\\
&\frac{S}{n}\frac{\tau}{\pi}c_n^{-\frac{m-1}{m}}\left(\frac{1}{2m-1}-\frac{1}{m-1}\right)=-c_n R(1+o(1)).
\end{align*}
Finally, for the solution of \eqref{24}, we can write 
\begin{align}
\label{25}
&c_n^\ast=\alpha_n^\ast (1+o(1)), \ \ \alpha_n^\ast=\left[\frac{S}{n R}\frac{\tau}{\pi}\frac{m}{(2m-1)(m-1)}\right]^\frac{m}{2m-1}.
\end{align}
Now, using the identity $(\beta\in\N,\,\beta>1)$
$$\sum_{l>N_n}\frac{1}{l^\beta}=\frac{1}{N_n^{\beta-1}}\int_1^{+\infty}\frac{1}{x^\beta}\d x\cdot (1+o(1)),\,n\longrightarrow+\infty,$$
for $\beta=2$
$$\sum_{l>N_n}\frac{1}{l^2}=\frac{1}{N_n}\cdot (1+o(1)),\,n\longrightarrow+\infty,$$
calculate
\begin{align*}
&H(c_n^\ast)=\frac{S}{n}\left[(c_n^\ast)^2\left(\frac{\pi}{\tau}\right)^{2(m-1)}\frac{N_n^{2m-1}}{2m-1}-\right.\left.2c_n^\ast\left(\frac{\pi}{\tau}\right)^{m-2}\frac{N_n^{m-1}}{m-1}-\left(\frac{\tau}{\pi}\right)^2\frac{1}{N_n}\right](1+o(1))\\
&+(c_n^\ast)^2R=\frac{S}{n}\frac{\tau}{\pi}\left[(c_n^\ast)^2\frac{(c_n^\ast)^{-\frac{2m-1}{m}}}{2m-1}-2c_n^\ast\frac{(c_n^\ast)^{-\frac{m-1}{m}}}{m-1}-(c_n^\ast)^\frac{1}{m}\right](1+o(1))+(c_n^\ast)^2R=\\
&=\frac{S}{n}\frac{\tau}{\pi}(c_n^\ast)^{\frac{1}{m}}\frac{-2m^2}{(2m-1)(m-1)}(1+o(1))+(c_n^\ast)^2R=\\
&=(-2m)R (c_n^\ast)^{\frac{1}{m}}(c_n^\ast)^{\frac{2m-1}{m}}(1+o(1))+(c_n^\ast)^2R=\\
&=-(2m-1)(\alpha_n^\ast)^2 R(1+o(1)),
\end{align*}
where we have used the relation \eqref{25}. Now, choosing the sequence $c_n=\alpha_n^\ast$ for the definition of the estimator in \eqref{29}, we obtain from \eqref{28} 
\begin{align}
\label{26}
&\sup_{\Lambda\in\F_m^{per}(R,S)}\left(\Ex_\Lambda\|\tilde\Lambda_n-\Lambda\|^2-\Ex_\Lambda\|\hat\Lambda_n-\Lambda\|^2\right)\leq\nonumber\\
&\leq-(2m-1)(\alpha_n^\ast)^2 R(1+o(1))+\frac{{\rm C}}{n}\max_l\frac{|1-|\tilde K_{2l,n}|^2|}{\left(\frac{2\pi l}{\tau}\right)^m}.
\end{align}
If we show that
\begin{align}
\label{27}
\frac{1}{n}\max_l\frac{|1-\tilde K_{l,n}^2|}{\left(\frac{\pi l}{\tau}\right)^m}=o(n^{-\frac{2m}{2m-1}}),
\end{align}
then, since (see \eqref{19})
$$\Pi=(2m-1)(\alpha_n^\ast)^2 Rn^{\frac{2m}{2m-1}},$$
we get from \eqref{26}
\begin{align*}
\varlimsup_{n\rightarrow+\infty}n^{\frac{2m}{2m-1}}\sup_{\Lambda\in\F_m^{per}(R,S)}\left(\Ex_\Lambda\|\tilde\Lambda_n-\Lambda\|^2-\Ex_\Lambda\|\hat\Lambda_n-\Lambda\|^2\right)\leq-\Pi.
\end{align*}
This combined with the proposition will end the proof. To prove \eqref{27} recall that
\begin{align*}
&\tilde K_{l,n}=\left(1-\left|\frac{\pi l}{\tau}\right|^m \alpha_n^\ast\right)_+,\ \ \alpha_n^\ast=\left[\frac{S}{n R}\frac{\tau}{\pi}\frac{m}{(2m-1)(m-1)}\right]^\frac{m}{2m-1}.
\end{align*}
Therefore, for $m>1$ we have
\begin{align*}
\frac{1}{n}\max_l\frac{|1-\tilde K_{l,n}^2|}{\left(\frac{\pi l}{\tau}\right)^m}\leq\frac{2}{n}\max_l&\frac{1-\tilde K_{l,n}}{\left(\frac{\pi l}{\tau}\right)^m}=\frac{2}{n}\alpha_n^\ast=\frac{{\rm C}}{n^{\frac{3m-1}{2m-1}}}=o(n^{-\frac{2m}{2m-1}}).
\end{align*}

\textbf{Acknowledgements} The author is grateful to Y.A. Kutoyants for his suggestions and interesting discussions, also to A.S. Dalalyan and Y.K. Golubev for their fruitful comments.

\end{document}